\documentclass[12pt]{amsart}
\usepackage[]{amsmath, amsthm, amsfonts, amssymb, epsfig}
\newcommand\C{{\mathbb C}}

\renewcommand\O{\Omega}
\newcommand\Obar{\overline{\Omega}}

\numberwithin{equation}{section}

\begin{document}

\title[Two-connected domains]
{A Riemann mapping theorem for two-connected domains in the plane}
\author[Bell, Deger, Tegtmeyer]
{Steven R. Bell, Ersin Deger, \& Thomas Tegtmeyer}

\address[]{Mathematics Department, Purdue University, West Lafayette,
IN  47907}
\email{bell@math.purdue.edu}

\address[]{Mathematics Department, University of Cincinnati, Cincinnati,
OH  45221}
\email{ersin.deger@uc.edu}

\address[]{Mathematics Department, Truman State University,
Kirksville, MO 63501}
\email{tegt@truman.edu}


\begin{abstract}
We show how to express a conformal map $\Phi$ of a general two connected
domain in the plane such that neither boundary component is a point to a
representative domain of the form ${\mathcal A}_r=\{z: |z+1/z|<2r\}$, where
$r>1$ is a constant.  The domain ${\mathcal A}_r$ has the virtue of having
an explicit algebraic Bergman kernel function, and we shall explain why
it is the best analogue of the unit disc in the two connected setting.
The map $\Phi$ will be given as a simple and explicit algebraic function
of an Ahlfors map of the domain associated to a specially chosen point.
It will follow that the conformal map $\Phi$ can be found by solving
the same extremal problem that determines a Riemann map in the simply
connected case.
\end{abstract}

\maketitle

\theoremstyle{plain}

\newtheorem {thm}{Theorem}[section]
\newtheorem {lem}[thm]{Lemma}

\hyphenation{bi-hol-o-mor-phic}
\hyphenation{hol-o-mor-phic}

\section{Main results}
\label{sec1}
The Riemann map $f_a$ associated to a point $a$ in a simply connected
domain $\O\ne\C$ in the complex plane maps $\O$ one-to-one onto the
unit disc $\{z: |z|<1\}$ and is the solution to an extremal problem: among
all holomorphic mappings of $\O$ into the unit disc, $f_a$ is the
map such that $f_a'(a)$ is real and as large as possible.  The
objects of potential theory associated to the unit disc are
particularly simple.  The Bergman and Szeg\H o kernel functions
are explicit and simple rational functions on the disc, for example.
The Poisson kernel is also very simple.  Consequently, the Bergman,
Poisson, and Szeg\H o kernels associated to $\O$ can be expressed very
simply in terms of a Riemann map.

It was noted in \cite{B4} that any two-connected domain $\O$ in the
plane such that neither boundary component is a point can be mapped
conformally via a mapping $\Phi$ to a representative domain of the form
$${\mathcal A}_r := \{z: |z+1/z|<2r\},$$
where $r$ is a constant greater than one.  See \cite{J-T} for a nice
proof of this fact which also proves a much more general theorem.
See Figure~1 for a picture of ${\mathcal A}_r$.  The dotted circle
in the figure is the unit circle.  Notice that $1/z$ is an automorphism
of ${\mathcal A}_r$ that fixes the unit circle.

In this paper, we show that the conformal map $\Phi:\Omega\to{\mathcal A}_r$
can be expressed simply in terms of an Ahlfors map of $\O$.  The Bergman
kernel associated to a representative domain ${\mathcal A}_r$ is an
algebraic function (see \cite{B4,B5}).  In fact, the Bergman kernel has
recently been written down explicitly in \cite{D}.
The Szeg\H o kernel is also algebraic and the Poisson kernel can
be expressed in terms of rather simple functions (see \cite{B6}).  It is
proved in \cite{B3}
that neither the Bergman nor Szeg\H o kernel can be rational in a
two-connected domain.  Thus, the representative domain ${\mathcal A}_r$
can be thought of as perhaps the best analogue of the unit
disc in the two-connected setting.  We shall show that the conformal
map $\Phi$ can be determined by solving {\it the same\/} extremal problem
as the one that determines a Riemann map in the simply connected case, but
we must be careful to choose the base point properly.  In fact, we
show how to find a point $a$ in $\O$ so that $\Phi$ is given by
$cf_a+\sqrt{c^2f_a^2-1}$ where $f_a$ is the Ahlfors mapping associated to
$a$, and $c$ is a complex constant that we will determine.  We will
explain along the way how $cf_a+\sqrt{c^2f_a^2-1}$ can be understood to
be a single valued holomorphic function on $\O$.  The Ahlfors
map $f_a$ is the solution to the extremal problem: among all holomorphic
mappings of $\O$ into the unit disc, $f_a$ is the map such that $f_a'(a)$
is real and as large as possible.  See \cite[p.~49]{B2} for the basic
properties of the Ahlfors map.

In order to choose $a$, we will inadvertently discover a new way to
compute the modulus of the domain $\O$.

\begin{center}
\epsfig{file=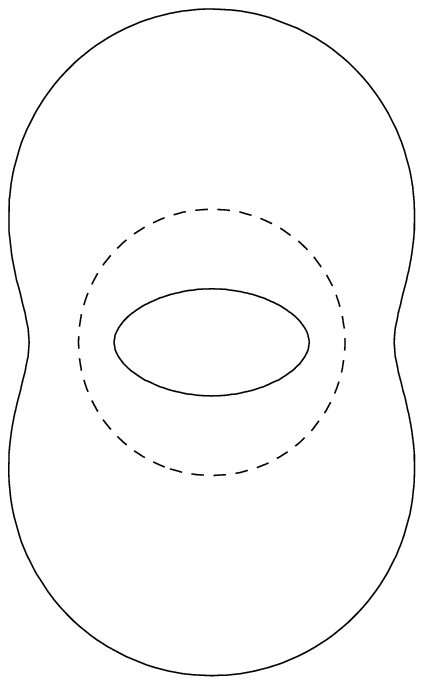, height=3in}
\end{center}

\begin{center}
Figure~1. The domain ${\mathcal A}_r$ with $r=1.05$.\\
\end{center}

Let $a_\rho$ denote the annulus $\{z:\frac{1}{\rho}<|z|<\rho\}$ when $\rho>1$, and
let
$$J(z):=\frac{1}{2}\left(z+\frac{1}{z}\right)$$
denote the Joukowsky map.  Note that $\frac{1}{r}J(z)$ is a proper holomorphic map of
${\mathcal A}_r$ onto the unit disc which is a $2$-to-one (counting multiplicities)
branched covering map.  Several authors have found explicit biholomorphic
maps between $a_\rho$ and ${\mathcal A}_r$ and have determined the
relationship between $\rho$ and $r$ (see \cite{J-T,JOT,C-M,D}).  Since we shall
need such a map, and since we can construct one in a minimum of space here using
the method from \cite{D}, we shall include the construction for completeness.
The biholomorphic mapping $\Psi: a_\rho \to {\mathcal A}_r$ will be an explicit
map involving the Jacobi {\it sn\/} function.  There is a standard conformal map
of the annulus $\{z: \frac{1}{\rho}<|z|<1\}$ onto the unit disc minus a slit from
$-L$ to $+L$ along the real axis.  See Nehari \cite[p.~293]{N}) for a construction
where the relationship between $L$ and $\rho$ is given explicitly.  The mapping is
given by a constant times
$$\text{\it sn}\left(K+\frac{2iK}{\pi}\log\rho z\, ;\, \rho^{-4}\right),$$
where $K$ is a constant that can be determined explicitly.  The map takes the
unit circle to the unit circle and the circle of radius~$1/\rho$ to the slit from
$-L$ to $L$.  By composing with an inversion of the annulus, we get a conformal
map that maps the annulus $\{z: \frac{1}{\rho}<|z|<1\}$ one-to-one onto the unit
disc minus a slit from $-L$ to $L$ that maps the unit circle to the slit from
$-L$ to $L$ and maps the circle of radius~$1/\rho$ to the unit circle.
The Joukowsky map $J(z)$ is a biholomorphic mapping of $\{z\in
{\mathcal A}_r : |z|>1\}$ onto the disc of radius $r$ minus a slit
from $-1$ to $+1$.  It maps the outer boundary of ${\mathcal A}_r$ to the unit
circle and the unit circle two-to-one onto the slit.  Thus, by scaling properly
and then composing, we can write down a biholomorphic mapping from
the annulus $\{z: \frac{1}{\rho}<|z|<1\}$ onto $\{z\in{\mathcal A}_r : |z|>1\}$
which maps the unit circle to the unit circle.  We may
now extend this map by reflection to a biholomorphic map from $a_\rho$ onto
${\mathcal A}_r$.  Finally, we may compose with the automorphism $1/z$ and a
rotation, if necessary, to obtain a conformal mapping
$\Psi: a_\rho\to{\mathcal A}_r$ such that the outer boundary of $a_\rho$ maps
to the outer boundary of ${\mathcal A}_r$ and such that $\Psi(i)=i$.  Note
that $\Psi$ is uniquely determined by these conditions.  We shall need to know
that $\Psi(-i)=-i$.  Since $-1/\Psi(-1/z)$ satisfies the same properties,
namely that it is a conformal map that takes the outer boundary to the outer
boundary and $\Psi(i)=i$, it must be that
$\Psi(z)\equiv -1/\Psi(-1/z)$.  From this, we deduce that
$\Psi(-i)\equiv -1/\Psi(-i)$, and hence that $\Psi(-i)$ is a square root
of $-1$.  Because $\Psi$ is one-to-one, it cannot be that $\Psi(-i)$ is
$i$, and so $\Psi(-i)=-i$.

It is well known that any two-connected domain $\O$ in the plane such that neither
boundary component is a point can be mapped biholomorphically via a map $F$ to a unique
annulus of the form $a_\rho$.  Let $C_1$ denote the unit circle.  We shall call
$\rho^2$ the modulus of $\O$ and we shall call the set $F^{-1}(C_1)$ the {\it median\/}
of the two-connected domain.  The median is the unique connected one-dimensional set
that is left invariant by every automorphism of the two-connected domain.  The
biholomorphic map $\Psi$ of $a_\rho$ onto ${\mathcal A}_r$ preserves the unit circle.
Hence, the median of ${\mathcal A}_r$ is also the unit circle.

Our main results will follow from two key lemmas, which we now list.

\begin{lem}
\label{lem1}
Suppose that $\O$ is a two-connected domain in the plane such that neither
boundary component is a point.  The Ahlfors map $f_a$ associated to any point $a$
in $\O$ is such that it has two distinct and simple branch points on the median
of $\O$.
\end{lem}

Lemma~\ref{lem1} was proved by McCullough and Mair in \cite{M-M} and
independently and by other means in \cite{Te}.

\begin{lem}
\label{lem2}
The Ahlfors map associated to the point $i$ in
${\mathcal A}_r$ is given by $\frac{1}{r}J(z)$.
\end{lem}

Lemma~\ref{lem2} was proved in \cite{D}.  For completeness, and because it
will not take much space, we will give alternate proofs of the lemmas in
the next section.

Suppose that $\O$ is a
two-connected domain in the plane such that neither boundary component
is a point.  We now describe the procedure to determine the mapping $\Phi$,
assuming the truth of the lemmas.  Pick any point $P$ in $\O$.  Lemma~\ref{lem1}
guarantees the existence of a point $a$ on the median of $\O$ such that $f_P'(a)=0$.
(Note that the Ahlfors map is eminently numerically computable if the boundary
of $\O$ is sufficiently smooth (see \cite{K-S,K-T,Tr,B1}), and the zeroes of a
function with only two simple zeroes can easily be located numerically as well.)
Now the Ahlfors map $f_a$ associated to the point $a$ is such that there are
two distinct points $p_1$ and $p_2$ on the median of $\O$ such that $f_a'(p_i)=0$
for $i=1,2$.  The points $p_i$ are also different from $a$ because $f_a'(a)\ne0$.
We now claim that there
exists an $r>1$ and a biholomorphic mapping from $\O$ to ${\mathcal A}_r$ such
that $\Phi(a)=i$.  We may also arrange for the outer boundary of $\O$ to get mapped
to the outer boundary of ${\mathcal A}_r$.  To see this, we use the
biholomorphic mapping $F$ of $\O$ onto an annulus $a_\rho$.
By composing with $1/z$, if necessary, we may map the outer boundary of $\O$ to
the outer boundary of $a_\rho$.  Since $a$ is on the median, $F(a)$ is on the
unit circle, and we may compose with a rotation so that $F(a)=i$.  Now we
may compose with the biholomorphic map $\Psi$ constructed above of
$a_\rho$ onto ${\mathcal A}_r$, noting that $i$ is fixed under $\Psi$ and that
outer boundaries are preserved, to obtain $\Phi$.  (At the moment, we can
only say that the number $r$ exists.  We will explain how to find it
momentarily.)  Lemma~\ref{lem2} states that $\frac{1}{r}J(z)$ is the Ahlfors
map of ${\mathcal A}_r$ associated to the point $i$.  The extremal property of
Ahlfors maps makes it easy to see that there is a unimodular complex number
$\lambda$ such that
\begin{equation}
\label{eqn1}
\frac{1}{r}J(\Phi(z))=\lambda f_a(z).
\end{equation}
Since the branch points of $J(z)$ occur at $\pm 1$, this last identity reveals
that $\Phi(p_i)=\pm 1$ for $i=1,2$.  We may assume that $\Phi(p_1)=-1$ and
$\Phi(p_2)=1$.  Since $J(-1)=-1$ and $J(1)=1$, equation~(\ref{eqn1}) now gives
$-\frac{1}{r}=\lambda f_a(p_1)$ and $\frac{1}{r}=\lambda f_a(p_2)$.
Either one of these two equations determines both $\lambda$ and $r$.
We can now solve equation~(\ref{eqn1}) locally for $\Phi$ to obtain
$$\Phi(z)=J^{-1}(c f_a(z)),$$
where $c=r\lambda$ and, of course, $J^{-1}(w)=w+\sqrt{w^2-1}$.
This formula extends by analytic continuation from any point in $\O$
that is not a branch point of $f_a$.  The apparent singularities at
the branch points are removable because the construction aligned the
two algebraic singularities of $J^{-1}$ to occur at the two images
of the branch points under $cf_a$.  We shall explain how to compute
$\Phi$ from its boundary values in \S\ref{sec3}.
The modulus $\rho^2$ of $\O$ can now be gotten from $r$ via any of
the formulas relating $\rho$ and $r$ proved in \cite{J-T,JOT,C-M,D}.

In the last section of this paper, we explore how the mapping $\Phi$
can be used to study the Bergman kernel of two-connected domains.

\section{Proofs of the Lemmas}
\label{sec2}

To prove Lemmas~\ref{lem1} and~\ref{lem2}, we will make use of the simple
fact that if $G:\O_1\rightarrow\O_2$ is a biholomorphism between bounded
domains, the Ahlfors maps transform according to the following formula.  Let
$f_j(z;w)$ denote the Ahlfors map $f_w$ associated to the point $w$ in $\O_j$.
The extremal property of the Ahlfors maps can be exploited to
show that
\begin{equation}
\label{eqn2}
f_1(z;w)=\lambda f_2\left(G(z);G(w)\right),
\end{equation}
where $\lambda$ is the unimodular constant that makes the derivative of
$f_1(z;w)$ in the $z$ variable at $z=w$ real and positive. 

\noindent
{\bf Proof of Lemma~\ref{lem1}.}
The biholomorphic map $F$ from $\O$ to $a_\rho$ together with
equation~(\ref{eqn2}) reveals that it is sufficient to prove Lemma~\ref{lem1}
for an annulus $a_{\rho}$ with base point $a>0$.  The function $h(z)=e^{iz}$
maps the horizontal strip $H=\{z\,:\,-\log \rho<{\rm Im}\, z<\log \rho\}$
onto $a_{\rho}$, takes the real line onto $C_1,$ and is $2\pi$-periodic.
Consider $g(z)=f_a(h(z))$.  Notice that $g$ is $2\pi$-periodic, maps  $H$
onto the unit disk, and maps the boundary of $H$ onto the unit circle.
It can be extended via a sequence of reflections to a function on $\C$.
Because of these reflections, $g$ has a second period of $(4\log\rho)i,$
making $g$ an elliptic function.  Because $f_a$ is two-to-one, $g$ is
a second order elliptic function. 

Any second order elliptic function $\phi$ satisfies the equation
$\phi(w-z)=\phi(z),$ where $w$ is the sum of two noncongruent zeroes of
$\phi$ (see \cite{H}).  From \cite{T-T}, we know that if $a$ is real,
$f_a$ has zeroes at $a$ and $-1/a$.  (We will prove this fact below to make
this paper self contained.)  It now follows that $g$ will have zeroes at
$-i\log a$ and $i\log a+\pi.$ Thus, we can choose $w=\pi.$ The mapping
$z\mapsto \pi-z$ has $\pi/2$ as a fixed point, while the point $3\pi/2$ gets
mapped to the congruent point $-\pi/2$.  These two points will be branch
points of $g,$ and $i$ and $-i$ will be the branch points of  $f_a$,
lying on the median of $a_{\rho}.$

\noindent
{\bf Proof of Lemma~\ref{lem2}.}
The biholomorphism $\Psi$ (see \S\ref{sec1}) from $a_{\rho}$ onto
$\mathcal A_r$ mapping $i$ to $i$ maps $-i$ to $-i$.  Let $f_i$
denote the Ahlfors map of $a_{\rho}$ with base point $i$.  From
\cite{T-T}, we know that $f_i$ will have $i$ and $-i$ as its zeroes.
(We include a proof of this fact in the next paragraph.)
Equation~(\ref{eqn2}) tells us that the Ahlfors map of ${\mathcal A}_r$
with base point $i$ will have $i$ and $-i$ as its zeroes.   Because
this Ahlfors map has the same zeroes as $\frac{1}{r}J(z),$ and both
functions map the boundary of $\mathcal A_r$ to the unit circle, they
are the same up to multiplication by a unimodular constant.  Since
$J'(i)=1>0$, this constant must be $1$.  Lemma~\ref{lem2} is proved.

We shall now flesh out the proofs of the lemmas by giving a shorter
and simpler proof than the one given in \cite{T-T} that the
Ahlfors map $f_a$ associated to a point $a$ in $a_\rho$ has simple
zeroes at $a$ and $-1/\bar{a}$.  Thus $f_a$ has simple zeroes at
$a$ and $-1/a$ if $a$ is real, and $f_i$ has simple zeroes at $i$
and $-i$.  The Ahlfors map $f_a$ always has a simple zero at $a$.
Since the Ahlfors map is two-to-one (counting multiplicities), $f_a$
must have exactly one other simple zero distinct from $a$.
The Ahlfors map is given via
\begin{equation}
\label{eqnSL}
f_a(z)=S(z,a)/L(z,a)
\end{equation}
where $S(z,a)$ is the Szeg\H o kernel and $L(z,a)$ is the
Garabedian kernel.  For the basic properties of the Szeg\H o and
Garabedian kernels, see \cite[p.~49]{B2}.  We note here that
$S(z,a)=\overline{S(a,z)}$ and that $L(z,a)\ne0$ if $a\in\O$ and
$z\in\Obar$ with $z\ne a$.  Also, $L(z,a)$ has a single simple pole
at $a$ in the $z$ variable, and $L(z,a)=-L(a,z)$.  It follows from
equation~(\ref{eqnSL}) that the other zero of $f_a$ is the one and
only zero of $S(z,a)$ in the $z$ variable.  Let $Z(a)$ denote
this zero.  It is proved in \cite{B4} that $Z(a)$
is a proper antiholomorphic self-correspondence of the domain.  In
the two-connected case that we are in now, it is rather easy to see that
$Z(a)$ is an antiholomorphic function of $a$ that maps $a_\rho$ one-to-one
onto itself.  To see this, let $S'(z,a)$ denote the derivative of $S(z,a)$
in the $z$ variable and notice that the residue theorem yields
$$Z(a)=\frac{1}{2\pi i}\int_{z\in b\O}\frac{ z S'(z,a)}{S(z,a)}\ dz.$$
(See \cite{B2} for proofs that $S(z,a)$ is smooth up to the boundary
and non-vanishing on the boundary in the $z$ variable when $a$ is held
fixed in $a_\rho$.)  It can be read off from this last formula that
$Z(a)$ is antiholomorphic in $a$.  Now hold $z$ fixed in $a_\rho$ and
let $a$ tend to a boundary point $a_b$.  Equation~(\ref{eqnSL})
shows that $f_a(z)$ tends to $S(z,a_b)/L(z,a_b)$.  But
$S(z,a_b)=-\frac{1}{i}L(z,a_b)T(a_b)$ where $T(a_b)$ is the complex
unit tangent vector at $a_b$ (see \cite[p.~107]{B2}), and
so $f_a(z)$ tends to $iT(a_b)$, a unimodular constant.  This shows
that the zeroes of $f_a(z)$ must tend to the boundary as $a$ tends to
the boundary.  Consequently, $Z(a)$ is a proper antiholomorphic self mapping
of the annulus.  The only such maps are one-to-one and onto (see \cite{M-R}).
We shall say that $Z(a)$ is an antiholomorphic automorphism of
$a_{\rho}$.  Note that $Z(a)$ must satisfy $Z(Z(a))=a$.  The only such
antiholomorphic automorphisms are $Z(a)=\lambda\bar{a}$, where $\lambda$
is a unimodular constant, and $Z(a)=1/\bar{a},$ and $Z(a)=-1/\bar{a}$.
Since $f_a$ does not have a double zero at $a$, $Z(a)$ cannot have a fixed
point.  Only $Z(a)=-1/\bar{a}$ is without fixed points in $a_{\rho},$ so
$f_a$ has zeroes at $a$ and $-1/\bar{a}$.  The proof is complete.

\section{How to compute $\Phi$ and the median of $\O$.}
\label{sec3}
Suppose that $\O$ is a bounded two connected domain bounded by two
non-intersecting $C^\infty$ Jordan curves.  The Kerzman-Stein-Trummer
method \cite{K-S,K-T} can be used to compute the boundary values of
Ahlfors mappings associated to $\O$.  (See \cite{B1} and \cite{B2},
Chapter~26, for a description of how to do this in the multiply
connected setting.)  Thus, we may pick a point $P$ in $\O$ and
compute $f_P$.  We shall explain how to compute the two branch
points of $f_P$ on the median of $\O$ below.  Hence, we may find
a point $a$ on the median of $\O$.  Let $f$ denote the Ahlfors
mapping associated to the point $a$ and let $p_1$
and $p_2$ denote the two simple zeroes of $f'$ on the median of $\O$.
We shall now show how to compute the points $p_1$ and $p_2$
and the function $\Phi$ using only the boundary values of the Ahlfors
map $f$.  The residue theorem yields that
$$p_1^n+p_2^n=\frac{1}{2\pi i}\int_{b\O} \frac{f''(z)z^n}{f'(z)}\ dz.$$
The two points can be determined from the values of
$p_1^n+p_2^n$ for $n$ equal to one and two.  Indeed the coefficients
$A$ and $B$ of the polynomial $z^2-Az+B=(z-p_1)(z-p_2)$ are given
as $A=p_1+p_2$ and $B=\frac{1}{2}[(p_1+p_2)^2-(p_1^2+p_2^2)]$, and
so $p_1$ and $p_2$ can be found via the quadratic formula.  Therefore,
the point $a$ itself and the branch points of $f_a$ can readily be found
numerically.

We may calculate $c=r\lambda=-1/f(p1)$ by the procedure
described in \S\ref{sec1}.  There are two holomorphic branches of
$J^{-1}(w)=w+\sqrt{w^2-1}$ near the boundary of the disc of radius~$r$
corresponding to the two choices of the square root.  One branch $j_o$
maps the circle of radius~$r$ one-to-one onto the outer boundary of ${\mathcal A}_r$
and the other branch $j_i$ maps the circle of radius~$r$ one-to-one onto the inner
boundary of ${\mathcal A}_r$.  The mapping $\Phi$ that maps $\O$ one-to-one
onto ${\mathcal A}_r$ can be determined from its boundary values on $\O$.
The boundary values of $\Phi$ on the outer boundary of $\O$ can be taken to
be $j_o(cf(z))$ and the boundary values on the inner boundary are then
$j_i(cf(z))$.  The extension of $\Phi$ to $\O$ can be gotten from the
boundary values of $\Phi$ via the Cauchy integral formula.

We now turn to a method to determine the median of $\O$.  Equation~(\ref{eqn1})
shows that $J(\Phi(z))= c f(z)$ where $c=\lambda r$ and $f$ is shorthand for
$f_a$.  This shows that $f$ maps the median of $\O$ to the line segment
${\mathcal L}$ from $-1/c$ to $1/c$.  Hence, the median of $\O$ is given
by $f^{-1}({\mathcal L})$.  We already know the points $p_1$, $p_2$, and $a$
on the median.  To compute other points $w$ in this set such that $f(w)=\tau$
for some $\tau$ in ${\mathcal L}$ where $\tau$ is not one of the endpoints
of ${\mathcal L}$, we may again use the residue theorem to obtain
$$w_1^n +w_2^n=\frac{1}{2\pi i}\int_{z\in b\O}\frac{f'(z)z^n}{f(z)-\tau}\ dz,$$
where $w_1$ and $w_2$ are the two points in $f^{-1}(\tau)$.  Finally,
we may use the Newton indentity for $n=2$ and the quadratic formula
as we did above to determine $p_1$ and $p_2$ to get $w_1$ and $w_2$.
Note that $p_1$ and $p_2$ correspond to the inverse images of the two
endpoints of ${\mathcal L}$, and so we may generate the rest of the median
by letting $\tau$ move between the endpoints.

\section{Other consequences.}
We describe in this section how the Bergman kernel of a two-connected domain
can be expressed rather concretely in terms of an Ahlfors map.
It was proved in \cite{D} that the Bergman kernel associated to ${\mathcal A}_r$
is given by
$$K_r(z,w)=C_1\frac{2k^2S(z,\bar w)+kC(z,\bar w)D(z,\bar w)+ C_2}
{z\bar w\sqrt{1-k^2J(z)^2}\sqrt{1-k^2J(\bar w)^2}},$$
where $k=1/r^2$, and $C_1$ and $C_2$ are constants that only depend on $r$, and
$$S(z,w)=-\left(\frac
{J(z)\frac{w^2-1}{2w}\sqrt{1-k^2J(z)^2}+
J(w)\frac{z^2-1}{2z}\sqrt{1-k^2J(w)^2}}
{1-k^2J(z)^2J(w)^2}\right)^2,$$
and
$$C(z,w)=
\frac
{-\frac{z^2-1}{2z}\frac{w^2-1}{2w}-J(z)J(w)\sqrt{1-k^2J(z)^2}\sqrt{1-k^2J(w)^2}}
{1-k^2J(z)^2J(w)^2},$$
and
$$D(z,w)=
\frac
{\sqrt{1-k^2J(z)^2}\sqrt{1-k^2J(w)^2}+k^2J(z)J(w)\frac{z^2-1}{2z}\frac{w^2-1}{2w}}
{1-k^2J(z)^2J(w)^2}.$$
We showed that the biholomorphic map $\Phi$ that we constructed above from a
two-connected domain to its representative domain ${\mathcal A}_r$ is such that
$\frac{1}{r}J(\Phi(z))=\lambda f_a(z)$, where $f_a$ is an Ahlfors map associated
to $\O$ for a point $a$ on the median of $\O$ and $\lambda$ is a unimodular
constant.  Differentiate this identity to see that
$$\Phi'(z)=2c f_a'(z)/(1-\Phi(z)^{-2})=
2c f_a'(z)/(1-J^{-1}(cf_a(z))^{-2}),$$
where $J^{-1}(w)=w+\sqrt{w^2-1}$.
The Bergman kernel associated to $\O$ is given by
$$\Phi'(z)K_r(\Phi(z),\Phi(w))\overline{\Phi'(w)}.$$
We may read off from this and the formulas above that the Bergman kernel
associated to $\O$ is $f_a'(z)\,\overline{f_a'(w)}$ times a rational function
of $f_a(z)$, and $\overline{f_a(w)}$, and simple algebraic functions of
$f_a$ of the form
$\sqrt{A f_a(z)^2-1}$ and the conjugate
of $\sqrt{Af_a(w)^2-1}$, where $A$ is equal to $c$, $\lambda$, or one.
The degree of the rational function does not depend on the modulus of
the domain, and the coefficients of the rational function depend in a
straightforward way on the modulus.

It was proved in \cite{B5} that the Bergman kernel associated to a
two connected domain such that neither boundary component is a point
can be expressed in terms of two Ahlfors maps $f_a$ and $f_b$ as
$$f_a'(z)R(z,w)\overline{f_a'(w)}$$
where $R(z,w)$ is a rational function of $f_a(z)$, $f_b(z)$,
$\overline{f_a(w)}$, and $\overline{f_b(w)}$.  The 
rational function has been proved to exist, but it has never been
given explicitly.  We are now in a position to write down this rational
function for the domain ${\mathcal A}_r$ using the points $a=i$ and
$b=1$.  It was proved in \cite{D} that the Ahlfors map associated
to the point $1$ of ${\mathcal A}_r$ is
$$g(z)=\frac{1}{2r}\frac{z-\frac{1}{z}}{\sqrt{1-k^2J(z)^2}}$$
where $k=1/r^2$.  The square root in the formula can be taken to be
the principal branch because the modulus of $kJ(z)$ is less than one
in ${\mathcal A}_r$.  We shall prove this result here because all
the ingredients are on the table.  It is a pleasant exercise to check
that $g$ is a two-to-one branched covering map of ${\mathcal A}_r$ onto
the unit disc.  Notice that $g$ has zeroes at $1$ and $-1$.
We proved above that the Ahlfors map $f_1$ has zeroes
at $1$ and $-1$.  Hence $f_1/g$ is a unimodular constant.   Checking
derivatives at the point $1$ yields that the unimodular constant is
$1$.  So $g$ is equal to the Ahlfors map $f_1$.   Recall that the
Ahlfors map $f_i$ associated to $i$ is just $\frac{1}{r}J(z)$.

To determine the rational function, we make note that the terms
$\frac{z^2-1}{2z}$ can be multiplied by unity in the form of
$$\frac{\sqrt{1-k^2J(z)^2}}{\sqrt{1-k^2J(z)^2}}$$
to see that
$$\frac{z^2-1}{2z}= r\frac{1}{2r}\left(z-\frac{1}{z}\right)=
rf_1(z)\sqrt{1-k^2J(z)^2}
=rf_1(z)\sqrt{1-r^2f_i(z)^2}.$$
The terms in the denominator in the expression for $K(z,w)$ can
be manipulated by multiplying by unity in the form of $f_i'/f_i'$,
noting that $f_i'(z)=\frac{1}{r}(1-\frac{1}{z^2})$,
to see that
$$z\sqrt{1-k^2J(z)}=\frac{2}{f_i'(z)} f_1(z)(1-r^2f_i(z)^2).$$
When these procedures are carried out in each of the expressions
that comprise the Bergman kernel, and when all the terms are
combined, the square roots all become squared and we obtain
$$K(z,w)=f_i'(z)
Q(f_i(z),f_1(z),\overline{f_i(w)},\overline{f_1(w)})
\overline{f_i'(w)},$$
where $Q(z_1,z_2,w_1,w_2)$ is a rational function of four
variables given by
$$C_1\left(\frac{\sigma(z_1,z_2,w_1,w_2)+\delta(z_1,z_2,w_1,w_2) +C_2}
{q(z_1,z_2,w_1,w_2)}\right),$$
where
$$\sigma(z_1,z_2,w_1,w_2)=-2\left(\frac{z_1w_2(1-r^2w_1^2)+w_1z_2(1-r^2z_1^2)}
{1-z_1^2w_1^2}
\right)^2$$
and
$$\delta(z_1,z_2,w_1,w_2)=-\left(\frac{(z_2w_2+z_1w_1)(1-r^2z_1^2)(1-r^2w_1^2)(1+z_1w_1z_2w_2) + C_2}
{(1-z_1^2w_1^2)^2}
\right),$$
and
$$q(z_1,z_2,w_1,w_2)=z_2w_2(1-r^2z_1^2)(1-r^2w_1^2).$$
That is more than twice as complicated as any formula on the
one-connected unit disc, but it could be worse.  A fascinating feature
that can be read off is that the coefficients depend very simply on $r$
and the degree of the rational function does not depend on $r$.

For the general two connected domain $\O$, we use the mapping $\Phi$
and the transformation formula for the Bergman kernels to read off that
$$K(z,w)=f_a'(z)R(z,w)\overline{f_a'(w)},$$
where $R(z,w)$ is a rational function of $f_a(z)$, $f_b(z)$,
$\overline{f_a(w)}$, and $\overline{f_b(w)}$.  We take
$a$ to be the point we chose in the construction of $\Phi$ and
$b$ to be equal to $p_2$.  We obtain a formula for the rational function
and we note that the degree does not depend on $r$.

It will be interesting to investigate if similar properties hold for
domains of higher connectivity using techniques of Crowdy and Marshall
\cite{C-M}.

\end{document}